\documentclass[11pt]{article}
\usepackage{amsfonts,eucal,amsmath,bm,color}
\usepackage{a4wide,hyperref}
\usepackage{amsfonts,amssymb,eucal,amsmath}

\newcommand{\nz}{\smallsetminus\{0\}}

\newtheorem{theorem}{Theorem}

\newtheorem{TheoremB}{Theorem}

\newtheorem{corollary}{Corollary}

\newtheorem{lemma}{Lemma}

\def\A{\mathbb{A}}

\def\R{\mathbb{R}}
\def\C{\mathbb{C}}

\def\Q{\mathbb{Q}}
\def\Z{\mathbb{Z}}
\def\N{\mathbb{N}}

\newcommand{\Rp}{\R^{+}}    

\def\cH{\mathcal{H}}

\def\cJ{\mathcal{J}}

\def\cR{\mathcal{R}}

\def\cP{\mathcal{P}}
\newcommand{\ra}{r_{\alpha}}

\newcommand{\ba}{\beta_\alpha}

\newcommand{\ve}{\varepsilon}

\renewcommand{\le}{\leqslant}
\renewcommand{\ge}{\geqslant}
\newcommand{\vv}[1]{{\mathbf{#1}}}

\newenvironment{proof}{\noindent\textsc{{Proof}}.}{\hfill\raisebox{-1ex}{$\boxtimes$}}

\newcommand{\Vol}{\operatorname{Vol}}

\newcounter{remark}
\newcounter{example}
\newcounter{question}
\newcounter{problem}
\newcounter{conjecture}
\begin{document}
\newlength{\mylength}
\newcommand{\sui}{[-\tfrac12,\tfrac12]}

\large

\title{{The distribution of close conjugate algebraic numbers}}

\author{Victor Beresnevich\footnote{EPSRC Advanced Research Fellow, grant no.EP/C54076X/1}\\ {\small\sc (York)} \and Vasili Bernik\\ {\small\sc (Minsk)} \and Friedrich G\"otze\\ {\small\sc (Bielefeld)} }

\date{}

\maketitle

{\footnotesize

\begin{abstract}
We investigate the distribution of real algebraic numbers of a fixed degree having a close conjugate number, the distance between the conjugate numbers being given as a function of their height. The main result establishes the ubiquity of such algebraic numbers in the real line and implies a sharp quantitative bound on their number. Although the main result is rather general it implies new estimates on the least possible distance between conjugate algebraic numbers, which improve recent bounds of Bugeaud and Mignotte. So far the results \`a la Bugeaud and Mignotte relied on finding explicit families of polynomials with clusters of roots. Here we suggest a different approach in which irreducible polynomials are implicitly tailored so that their derivatives assume certain values. The applications of our main theorem considered in this paper include generalisations of a theorem of Baker and Schmidt and a theorem of Bernik, Kleinbock and Margulis in the metric theory of Diophantine approximation.
\end{abstract}



}

\section{Introduction}\label{itro}

\subsection{Separation of conjugate algebraic numbers}

The question \emph{``How close to each other can two conjugate algebraic numbers of degree $n$ be?''} crops up in a variety of problems in Number Theory and in some applications. Over the past 50 years or so there has been found a number of upper and lower bounds for such distance. However, the exact answers are known in the case of degree $2$ and $3$ only. In order to set the sense in our discussion we now introduce some quantities.

Throughout this paper we deal with algebraic numbers in $\C$, the set of complex numbers. Let $n\ge 2$. Recall that complex algebraic numbers are called \emph{conjugate} (over $\Q$) if they are roots of the same irreducible (over $\Q$) polynomial with rational integer coefficients. Define $\kappa_n$ (respectively $\kappa^*_n$) to be the infimum of $\kappa$ such that the inequality
$$
 |\alpha_1-\alpha_2| > H(\alpha_1)^{-\kappa}
$$
holds for arbitrary conjugate algebraic numbers (respectively algebraic integers) $\alpha_1\not=\alpha_2$ of degree $n$ with sufficiently large height $H(\alpha_1)$. Here and elsewhere $H(\alpha)$ denotes the height of an algebraic number $\alpha$, which is the absolute height of the minimal polynomial of $\alpha$ over $\Z$.
Clearly, $\kappa^*_n\le\kappa_n$ for all $n$.

In 1964 Mahler \cite{Mahler-64:MR0166188} proved the upper bound $\kappa_n\le n-1$, which is apparently the best estimate up to date. It is an easy exercise to show that $\kappa_2=1$ (see, e.g. \cite{Bugeaud-Mignotte-09}). Furthermore, Evertse \cite{Evertse-04:MR2107951} proved that $\kappa_3=2$. In the case of algebraic integers $\kappa^*_2=0$ and $\kappa^*_3\ge 3/2$. The latter has been proved by Bugeaud and Mignotte \cite{Bugeaud-Mignotte-09} who have also shown that the equality $\kappa^*_3=3/2$ is equivalent to Hall's conjecture on the difference between integers $x^3$ and $y^2$. The latter is known to be the special case of the $abc$-conjecture of Masser and Oesterl\'e -- see \cite{Bugeaud-Mignotte-09} for further details and references.

For $n>3$ estimates for $\kappa_n$ are less satisfactory. At first Mignotte \cite{Mignotte-83:MR728976} showed that $\kappa_n,\kappa^*_n\ge n/4$ for all $n\ge3$. Recently Bugeaud and Mignotte \cite{Bugeaud-Mignotte-09, Bugeaud-Mignotte-2004:MR2096618} have shown that
$$
\begin{array}{lcl}
  \kappa_n\ge n/2 && \text{when \ $n\ge 4$ \ is even,} \\[0.5ex]
  \kappa^*_n\ge (n-1)/2 && \text{when \ $n\ge 4$ \ is even,} \\[0.5ex]
  \kappa_n\ge (n+2)/4 && \text{when \ $n\ge 5$ \ is odd,} \\[0.5ex]
  \kappa^*_n\ge (n+2)/4 && \text{when \ $n\ge 5$ \ is odd.}
\end{array}
$$

The above results are obtained by presenting explicit families of irreducible polynomials of degree $n$ whose roots are close enough. Bugeaud and Mignotte \cite{Bugeaud-Mignotte-09} point out that ``at present there is no general theory for constructing integer polynomials of degree at least four with two roots close to each other''. In this paper we shall make an attempt to address this issue. One particular consequence of our results is the following theorem that improves the lower bounds of Bugeaud and Mignotte in the apparently more difficult case of odd $n$:

\bigskip

\begin{theorem}\label{t:01}
For any $n\ge 2$ we have that $\min\{\kappa_n,\kappa^*_{n+1}\} \ge (n+1)/3$.
\end{theorem}

\bigskip

Theorem~\ref{t:01} will follow from a more general counting result -- Corollary~\ref{c:02} below. In fact, a lot more is established. We show that algebraic numbers of degree $n$ (algebraic integers of degree $n+1$) with a close conjugate form a `highly dense' (ubiquitous) subset in the real line, see Theorem~\ref{t:02}.

\subsection{The distribution of close conjugate algebraic numbers}

First some notation. Throughout, $\#S$ stands for the cardinality of $S$ and $\lambda$ will denote Lebesgue measure in $\R$. Given an interval $J\subset\R$, $|J|$ will denote the length of $J$. Also, $B(x,\rho)$ will denote the interval in $\R$ centred at $x$ of radius $\rho$. By $\ll$ ($\gg$) we will mean the Vinogradov symbols with implicit constant depending on $n$ only. We shall write $a\asymp b$ when the inequalities $a\ll b$ and $a\gg b$ hold simultaneously.

Let $n\ge 2$ be an integer, $\mu\ge0$, $0<\nu<1$ and $Q>1$. Let $\A_{n,\nu}(Q,\mu)$ be the set of algebraic numbers  $\alpha_1\in\R$ of degree $n$ and height $H(\alpha_1)$ satisfying
\begin{equation}\label{ne:000}
\nu Q\le H(\alpha_1)\le \nu^{-1}Q
\end{equation}
and
\begin{equation}\label{ne:001}
 \nu\,Q^{-\mu}\le |\alpha_1-\alpha_2|\le \nu^{-1} Q^{-\mu}\qquad\text{for some $\alpha_2\in\R$, conjugate to $\alpha_1$}.
\end{equation}
Similarly we define $\A^*_{n,\nu}(Q,\mu)$ to be the set of algebraic integers $\alpha_1\in\R$ of degree $n+1$ and height $H(\alpha_1)$ satisfying (\ref{ne:000}) and (\ref{ne:001}).
Before we state our main result let us agree that $\A^\circ_{n,\nu}(Q,\mu)$ will refer to any of the sets $\A_{n,\nu}(Q,\mu)$ and $\A^*_{n,\nu}(Q,\mu)$.

\bigskip

\begin{theorem}\label{t:02}
For any $n\ge 2$ there is a constant
$\nu>0$ depending on $n$ only with the following
property. For any $\mu$ satisfying
\begin{equation}\label{e:007}
0<\mu\le \frac{n+1}{3}
\end{equation}
and any interval $J\subset [-\tfrac12,\tfrac12]$, for all sufficiently large $Q$
\begin{equation}\label{vb1}
    \lambda\left(\bigcup_{\alpha_1\in\A^\circ_{n,\nu}(Q,\mu)}B(\alpha_1,Q^{-n-1+2\mu})\cap J\right)\ge \tfrac34|J|.
\end{equation}
\end{theorem}

\bigskip

\begin{remark}
The constant $\tfrac34$ in the right hand side of (\ref{vb1}) is not `critical' and can be replaced by any positive number $<1$.
\end{remark}

\begin{remark}
In principle, the above theorem holds in the case $\mu=0$, though it is more delicate to ensure that $\alpha_2$ is real. Also note that in the case $\mu=0$ the optimal distribution of real algebraic numbers of degree $n$ (algebraic integers of degree $n+1$) was first established in \cite{Beresnevich-99:MR1709049} (respectively in \cite{Bugeaud-?2}). The above result is stated for the unit symmetric interval $\sui$. However, using shifts by an integer it can be extended to an arbitrary interval in $\R$ -- see \cite{Beresnevich-99:MR1709049} for appropriate technique.
\end{remark}

\bigskip

\begin{corollary}\label{c:01}
For any $n\ge 2$ there is a positive constant
$\nu$ depending on $n$ only such that for any $\mu$ satisfying $(\ref{e:007})$ and any interval $J\subset [-\tfrac12,\tfrac12]$, for all sufficiently large $Q$
\begin{equation}\label{vb3}
\#\big(\A^\circ_{n,\nu}(Q,\mu)\cap J\big)\ge \tfrac12Q^{n+1-2\mu}|J|.
\end{equation}
\end{corollary}

\bigskip

\begin{proof}
Obviously if $B(\alpha_1,Q^{-n-1+2\mu})\cap \frac12J \not=\varnothing$ then $\alpha_1\in J$ provided that $Q$ is sufficiently large. Then, using (\ref{vb1}) we obtain
$$
\#\big(\A^\circ_{n,\nu}(Q,\mu)\cap J\big)2Q^{-n-1+2\mu}\ge
\lambda\left(\bigcup_{\alpha_1\in\A^\circ_{n,\nu}(Q,\mu)}\hspace*{-3ex}B(\alpha_1,Q^{-n-1+2\mu})\cap \tfrac12J\right)\stackrel{\eqref{vb1}}{\ge} \tfrac14|J|,
$$
whence (\ref{vb3}) readily follows.
\end{proof}

\bigskip

\begin{corollary}\label{c:02}
Let $n\ge 2$. Then for all sufficiently large $Q>1$ there are $\gg Q^{\frac{n+1}3}$ real algebraic numbers $\alpha_1$ of degree $n$ $($real algebraic integers $\alpha_1$ of degree $n+1)$ with height $H(\alpha_1)\asymp Q$ such that
\begin{equation}\label{e:009}
|\alpha_1-\alpha_2|\asymp Q^{-\frac{n+1}{3}}\qquad\text{for some $\alpha_2\in\R$, conjugate to $\alpha_1$.}
\end{equation}
\end{corollary}

\bigskip

Corollary~\ref{c:02} follows from Corollary~\ref{c:01} on taking $\mu$ to be $\frac{n+1}{3}$. As a consequence of Corollary~\ref{c:02} we obtain Theorem~\ref{t:01}.

\section{Auxiliary lemmas}

The following auxiliary statement established in \cite[Theorem~5.8]{Beresnevich-SDA1} is the crucial ingredient of the proof of all the results of this paper.

\bigskip

\begin{lemma}\label{l:01}
Let $f_0,\dots,f_n$ be real analytic linearly
independent over $\R$ functions defined on an interval $I\subset\R$. Let $x_0\in I$ be a point such that the Wronskian $W(f_0,\dots,f_n)(x_0)\not=0$. Then there is an interval $I_0\subset I$ centred at $x_0$ and positive constants $C$ and $\alpha$
satisfying the following property. For any interval $J\subset I_0$ there is a constant $\delta=\delta_J$ such that for any positive $\theta_0,\dots,\theta_{n}$
\begin{equation}\label{e:016}
\lambda\left\{ x\in J  :
\begin{array}{l}
    \exists\ (a_0,\dots,a_n)\in\Z^{n+1}\nz \ \text{satisfying} \\[0.2ex]
    |a_0f_0^{(i)}(x)+\dots+a_nf_n^{(i)}(x)|< \theta_i\ \forall\ i=\overline{0,n}
\end{array}
\right\} \le
 C\Big(1+\Big(\frac{\Theta}{\delta}\Big)^{\alpha}\Big)\theta^\alpha|J|,
\end{equation}
where
\begin{equation}\label{e:017}
\theta=(\theta_0\dots\theta_{n})^{1/(n+1)}\qquad\text{and}\qquad\Theta:=\max_{1\le
    r\le n}\frac{\theta_0\cdots\theta_{r-1}}{\theta^r}.
\end{equation}
\end{lemma}

\bigskip

With the view to the applications we have in mind we now
estimate $\Theta$.

\begin{lemma}\label{l:02}
Assume that $\theta\le1$ and
assume that for some index $m\le n$ we have
\begin{equation}\label{e:019}
\theta_0, \dots, \theta_{m-1} \le k \qquad\text{and}\qquad \theta_m,  \dots, \theta_n\ge k^{-1}\qquad\text{for some real $k\ge1$.}
\end{equation}
Then
\begin{equation}\label{e:020}
\Theta \le k^{n-1}\max\left\{\frac{\theta_0}{\theta_0\dots\theta_n},\frac{1}{\theta_n}\right\}\,.
\end{equation}
\end{lemma}

\bigskip

\begin{proof} By the assumption that $\theta\le1$, for
all $r\in\{1,\dots,n\}$ we have $\theta^r\ge
\theta^{n+1}=\theta_0\dots\theta_n$. Therefore $\Theta$ satisfies
\begin{equation}\label{e:021}
\Theta\le \frac{1}{\theta_0\dots\theta_n}\max_{1\le r\le n}\theta_0\cdots\theta_{r-1}.
\end{equation}
In view of (\ref{e:019}) it is readily seen that
\begin{equation}\label{vb+1}
\max_{1\le r\le m}\theta_0\cdots\theta_{r-1}\le k^{n-1}\max_{1\le r\le m} \theta_0\frac{\theta_1}{k}\cdots\frac{\theta_{r-1}}{k}\stackrel{\eqref{e:019}}{=} k^{n-1}\theta_0
\end{equation}
and
\begin{equation}\label{vb+2}
\max_{m< r\le n}\theta_0\cdots\theta_{r-1}\le\max_{m< r\le n}\prod_{i=0}^{m-1}\theta_i\prod_{i=m}^{r-1}k\theta_{i}
\stackrel{\eqref{e:019}}{=} \prod_{i=0}^{m-1}\theta_i\prod_{i=m}^{n-1}k\theta_{i}\le k^{n-1}\theta_0\cdots\theta_{n-1}.
\end{equation}
Combining (\ref{vb+1}) and (\ref{vb+2}) with (\ref{e:021}) gives (\ref{e:020}).
\end{proof}

\bigskip

We will be using Lemma~\ref{l:01} with $f_i(x)=x^i$
$(0\le i\le n)$. In this case the Wronskian $W(f_0,\dots,f_n)$ identically equals $n!$ and Lemma~\ref{l:01} is applicable to a neighborhood of any point $x_0\in\R$. The system of inequalities in (\ref{e:016}) becomes
\begin{equation}\label{e:018}
    |P(x)| < \theta_0,\quad
    |P'(x)| < \theta_1,\ \dots,\quad
    |P^{(n)}(x)| <  \theta_n,
\end{equation}
where $P(x)=a_0+a_1x+\dots+a_nx^n$ is a non-zero integral polynomial of degree at most $n$ and the set in the left hand side of (\ref{e:016}) is simply
\begin{equation}\label{e:018+}
A_n(J;\theta_0,\dots,\theta_{n}):=\left\{x\in J:
\begin{array}{l}
\text{(\ref{e:018}) holds for some}\\[0.3ex]
P\in\Z[x]\nz,\deg P\le n
\end{array}
\right\}\,.
\end{equation}
Then, combining Lemmas~\ref{l:01} and \ref{l:02} and using pretty standard compactness argument (e.g., \cite[proof of Lemma~6]{Beresnevich-Dickinson-Velani-07:MR2373145}) give

\bigskip

\begin{lemma}\label{l:03}
There are constants $C>0$ and $\alpha>0$ depending on $n$ only such that for any interval $J\subset \sui$ there is a constant $\delta_J>0$ such that for any positive numbers $\theta_0,\dots,\theta_{n}$ satisfying $\theta=(\theta_0\dots\theta_n)^{1/(n+1)}\le1$ and $(\ref{e:019})$ we have that
\begin{equation}\label{e:022}
\lambda\big(A_n(J;\theta_0,\dots,\theta_{n})\big)  \le
C\left(1+\frac{k^{\alpha(n-1)}}{\delta_J^\alpha}
\max\left\{\frac{\theta_0}{\theta_0\dots\theta_n},\frac{1}{\theta_n}\right\}^\alpha\right)\theta^\alpha|J|.
\end{equation}
\end{lemma}

\section{Tailored polynomials}

Let $\xi_0,\dots,\xi_{n}\in\Rp$ satisfy the conditions
\begin{equation}\label{e:039}
\begin{array}{l}
 \xi_i\ll 1\quad\text{when }0\le i\le m-1,\\[1ex]
 \xi_i\gg 1 \quad\text{when }m\le i\le n,\\[1ex]
 \xi_0<\ve,\qquad \xi_n>\ve^{-1}
\end{array}
\end{equation}
for some $0<m\le n$ and $\ve>0$, where the implied constants depend on $n$ only. Assume also that
\begin{equation}\label{e:040}
\prod_{i=0}^n\xi_i=1.
\end{equation}
The following lemma lies at the heart of the proof of Theorem~\ref{t:02}. It enable us to tailor irreducible polynomials which assume certain values of derivatives. Of course, there is a connection with Taylor's formula too. Hence, we call them \emph{tailored polynomials}.

\bigskip

\begin{lemma}\label{l:04}
For every $n\ge 2$ there are positive constants $\delta_0$ and $c_0$ depending on $n$ only with the following property. For any interval $J\subset\sui $ there is a sufficiently small $\ve=\ve(n,J)>0$ such that for any $\xi_0,\dots,\xi_{n}$ satisfying $(\ref{e:039})$ and $(\ref{e:040})$ there is a measurable set $G_J\subset J$ satisfying
\begin{equation}\label{e:024}
    \lambda(G_J)\ge\tfrac34|J|
\end{equation}
such that for every $x\in G_J$ there are $n+1$ linearly independent primitive irreducible polynomials $P\in\Z[x]$ of degree exactly $n$ such that
\begin{equation}\label{e:041}
 \delta_0 \xi_i \le |P^{(i)}(x)|  \le  c_0 \xi_i\qquad\text{for all }i=0,\dots,n\,.
\end{equation}
\end{lemma}

\bigskip

\begin{proof}
Let $n\ge 2$ and let $\xi_0,\dots,\xi_{n}$ be given and satisfy $(\ref{e:039})$ and $(\ref{e:040})$ for some $m$ and $\ve$. Let $J\subset\sui $ be any interval and $x\in J$. Consider the system of inequalities
\begin{equation}\label{e:026}
 |P(x)|  \le  \xi_i\qquad\text{when }0\le i\le n\,,
\end{equation}
where $P(x)=a_nx^n+\dots+a_1x+a_0$. Let $B_x$ be the set of $(a_0,\dots,a_n)\in\R^{n+1}$ satisfying (\ref{e:026}). Clearly, $B_x$ is a convex body in $\R^{n+1}$ symmetric about the origin. In view of (\ref{e:040}), the volume of this body equals $2^{n+1}\prod_{i=1}^ni!^{-1}$. Let
$\lambda_0\le\lambda_1\le\dots\le \lambda_n$ be the successive
minima of $B_x$. Clearly, $\lambda_i=\lambda_i(x)$ is a function of $x$. By Minkowski's theorem for successive minima,
$$
 \frac{2^{n+1}}{(n+1)!} \ \le \ \lambda_0\dots\lambda_n\Vol{B_x}\ \le \ 2^{n+1}.
$$
Substituting the value of $\Vol{B_x}$ gives
$\lambda_0\dots\lambda_n\ \le \ \prod_{i=1}^ni!$. Therefore, since
$\lambda_0\le\dots\le\lambda_n$, we get that
\begin{equation}\label{e:027}
\lambda_n\le \lambda_0^{-n}\prod_{i=1}^ni!.
\end{equation}

Our next goal is to show that $\lambda_0$ is bounded below by a
constant unless $x$ belongs to a small subset of $J$. Let
$E_\infty(J,\delta_1)$ be the set of $x\in J$ such that
$\lambda_0=\lambda_0(x)\le\delta_1$, where $\delta_1<1$. By the definition of
$\lambda_0$, there is a non-zero polynomial $P\in\Z[x]$, $\deg P\le
n$ satisfying
\begin{equation}\label{e:028}
 |P^{(i)}(x)|  \le \delta_1 \xi_i\qquad  (0\le i\le n).
\end{equation}
Let $\theta_0=\delta_1\xi_0$ and $\theta_i=\xi_i$ ($1\le i\le n$). Then $E_\infty(J,\delta_1) \subset A_n(J;\theta_0,\dots,\theta_{n})$ -- see (\ref{e:018+}) for the definition of $A_n(\cdot)$. In view of (\ref{e:039}) and (\ref{e:040}), Lemma~\ref{l:03} is applicable. For this choice of $\theta_0,\dots,\theta_n$ we have $\theta=\delta_1^{1/(n+1)}$. Then
$$
\lambda(E_\infty(J,\delta_1)) \le \lambda\big(A_n(J;\theta_0,\dots,\theta_{n})\big)   \ll \left(1+\frac{1}{\delta_J^\alpha}\max\left\{\frac{\delta_1\xi_0}{\delta_1},
\frac{1}{\xi_n}\right\}^{\alpha}\right)
\delta_1^{\frac{\alpha}{n+1}}|J|.
$$
By (\ref{e:039}), $\max\{\xi_0,\xi_n^{-1}\}<\ve$. Therefore $\mu(E_\infty(J,\delta_1))\ll \delta_1^{\alpha/(n+1)}|J|$ provided that $\ve<\delta_J$. Then there is a sufficiently small $\delta_1$ depending on $n$ only such that
\begin{equation}\label{e:029}
\lambda(E_\infty(J,\delta_1))\le \tfrac1{4n+8}|J|.
\end{equation}
By construction, for any $x\in
J\setminus E_\infty(J,\delta_1)$ we have that
\begin{equation}\label{e:030}
 \lambda_0\ge \delta_1.
\end{equation}
Combining (\ref{e:027}) and (\ref{e:030}) gives
\begin{equation}\label{e:031}
 \lambda_n\le c_1:=\delta_1^{-n}\prod_{i=1}^ni!,
\end{equation}
where $c_1$ depends on $n$  only. By the definition of $\lambda_n$, there are $(n+1)$ linearly
independent integer points $\vv a_j=(a_{0,j},\dots,a_{n,j})$ $(0\le
j\le n)$ lying in the body $\lambda_n B_x\subset c_1 B_x$. In other
words, the polynomials $P_j(x)=a_{n,j}x^n+\dots+a_{0,j}$ $(0\le j\le n)$ satisfy the system of inequalities
\begin{equation}\label{e:032}
 |P_j^{(i)}(x)|  \le  c_1 \xi_i\qquad (0\le i\le n).
\end{equation}

Let $A=(a_{i,j})_{0\le i,j\le n}$ be the integer matrix composed from the integer points $\vv a_j$ $(0\le j\le n)$. Since all these points are contained in the body $c_1B_x$, we have that $|\det A|\ll\Vol(B_x)\ll 1$. That is $|\det A| < c_2$ for some constant $c_2$ depending on $n$ only. By Bertrand's postulate, choose a prime number $p$ satisfying
\begin{equation}\label{e:033}
 c_2\le p \le 2c_2.
\end{equation}
Therefore, $|\det A|<p$. Since $\vv a_0,\dots,\vv a_n$ are linearly independent and integer, $|\det A|\ge1$. Therefore, $\det
A\not\equiv 0\pmod p$ and the following system
\begin{equation}\label{e:034}
A \overline t\equiv \overline b\pmod p
\end{equation}
has a unique non-zero integer solution $\overline t={}^t(t_0,\dots,t_n)\in[0,p-1]^{n+1}$, where $\overline  b:={}^t(0,\dots,0,1)$ and ${}^t$ denotes transposition. For $l=0,\dots,n$ define $\overline r_l={}^t(1,\dots,1,0,\dots,0)\in\Z^{n+1}$, where the number of zeros is $l$. Since $\det A\not\equiv 0\pmod p$, for every $l=0,\dots,n$ the following system
\begin{equation}\label{add}
A\overline \gamma\equiv -\frac{A\overline t-\overline b}{p}+\overline r_l\pmod{p}
\end{equation}
has a unique non-zero integer solution $\overline \gamma=\overline  \gamma_l\in[0,p-1]^{n+1}$. Define $\overline\eta=\eta_l:=\overline t+ p\overline\gamma_l$ \ ($0\le l\le n$). Consider the $(n+1)$ polynomials of the form
\begin{equation}\label{add2}
 P(x)=a_nx^n+\dots+a_0:=\sum_{i=0}^n\eta_iP_i(x)\in\Z[x],
\end{equation}
where $(\eta_0,\dots,\eta_n)=\eta$ which of course depends on the parameter $l\in\{0,\dots,n\}$. Since $\overline r_0,\dots,\overline r_n$ are linearly independent, it is easily seen that $\overline \eta_0,\dots,\overline \eta_n$ are linearly independent.
Hence the polynomials given by (\ref{add2}) are linearly independent and so are non-zero.

Observe that $A\overline\eta$ is actually the column ${}^t(a_0,\dots,a_n)$ of coefficients of $P$. By construction, $\overline \eta\equiv\overline t\pmod p$ and therefore $\overline \eta$ is also a solution of (\ref{e:034}). Then, since $\overline b={}^t(0,\dots,0,1)$ and $A\overline\eta\equiv\overline b\pmod p$,
we have that $a_n\not\equiv0\pmod p$ and $a_i\equiv0\pmod p$ for $i=0,\dots,n-1$. Furthermore, by (\ref{add}), we have that $A\overline\eta\equiv\overline b+p\overline r_l\pmod{p^2}$. Then, on substituting the values of $\overline b$ and $\overline r_l$ into this congruence one readily verifies that $a_0\equiv p\pmod{p^2}$ and so $a_0\not\equiv0\pmod{p^2}$. By Eisenstein's criterion, $P$ is irreducible.

Since both $\overline t$ and $\overline \gamma_l$ lie in $\in[0,p-1]^{n+1}$ and $\overline \eta=\overline t+p\overline \gamma_l$, it is readily seen that $|\eta_i|\le p^2$ for all $i$. Therefore, using (\ref{e:032}) and (\ref{e:033}) we obtain that
\begin{equation}\label{e:035}
 |P^{(i)}(x)|  \le  c_0 \xi_i\qquad (0\le i\le n)
\end{equation}
with $c_0=4(n+1)c_1c_2^2$. Without loss of generality we may assume that the $(n+1)$ linearly independent polynomials $P$ constructed above are primitive (that is the coefficients of $P$ are coprime) as otherwise the coefficients of $P$ can be divided by their greatest common multiple. Clearly such division would not affect the validity of (\ref{e:035}). Thus, $P\in\Z[x]$ are primitive irreducible polynomials of degree $n$ which satisfy the right hand side of (\ref{e:041}). The final part of the proof is aimed at establishing the left hand side of (\ref{e:041}). The arguments are applied to every of the polynomials $P$ we have constructed.

Let $\delta_0>0$ be a sufficiently small parameter depending on $n$. For every $j=\overline{0,n}$ let $E_j(J,\delta_0)$ be the set of $x\in
J$ such that there is a non-zero polynomial $R\in\Z[x]$, $\deg R\le
n$ satisfying
\begin{equation}\label{e:036}
 |R^{(i)}(x)|  \le  \delta_0^{\delta_{i,j}}c_0^{1-\delta_{i,j}}\xi_i,
\end{equation}
where $\delta_{i,j}$ equals $1$ if $i=j$ and $0$ otherwise.
Let $\theta_i=\delta_0^{\delta_{i,j}}c_0^{1-\delta_{i,j}}\xi_i$. Then $E_j(J,\delta_0) \subset A_n(J;\theta_0,\dots,\theta_{n})$. In view of (\ref{e:039}) and (\ref{e:040}), Lemma~\ref{l:03} is applicable provided that $\ve<\min\{c_0^{-1},c_0\delta_0\}$. Then, by Lemma~\ref{l:03},
$$
\lambda(E_j(J,\delta_0)) \ll
\left(1+\frac{1}{\delta_J^\alpha}\max\left\{\frac{c_0\xi_0}{c_0^n\delta_0},
\frac{1}{\delta_0c_0\xi_n}\right\}^{\alpha}\right)(\delta_0c_0^n)^{1/(n+1)}|J|\,.
$$
It is readily seen that the above maximum is $\le\delta_J$ if $\ve<\delta_J\delta_0c_0$. Then
 \begin{equation}\label{e:038}
\lambda(E_j(J,\delta_0))\le \tfrac1{4n+8}|J|
\end{equation}
provided that $\ve<\min\{\delta_J\delta_0c_0,c_0^{-1},c_0\delta_0\}$ and $\delta_0=\delta_0(n)$ is sufficiently small. By construction, for any $x$ in the set $G_J$ defined by
$$
G_J:=J\setminus\left(\bigcup_{j=0}^nE_j(J,\delta_0)\cup E_\infty(J,\delta_1)\right)
$$
we must necessarily have that $|P^{(i)}(x)|\ge\delta_0\xi_i$ for all $i=0,\dots,n$, where $P$ is the same as in (\ref{e:035}). Therefore, the left hand side of (\ref{e:041}) holds for all $i$.
Finally, observe that
$$
\lambda(G_J)\ge |J|-\sum_{i=0}^n\lambda(E_i(J,\delta_0))-\lambda(E_\infty(J,\delta_1))\!\stackrel{\eqref{e:029}\,\&\,\eqref{e:038}}{\ge}\!
|J|-(n+2)\,\tfrac1{(4n+8)}|J|=\tfrac34|J|.
$$
The latter verifies (\ref{e:024}) and completes the proof.
\end{proof}

\section{Tailored monic polynomials}

The following is the analogue of Lemma~\ref{l:04} for the case of monic polynomials.

\bigskip

\begin{lemma}\label{l:05}
For every $n\ge 2$ there are positive constants $\delta_0$ and $c_0$ depending on $n$ only with the following property. For any interval $J\subset\sui $ there is a sufficiently small $\ve=\ve(n,J)>0$ such that for any positive $\xi_0,\dots,\xi_{n}$ satisfying $(\ref{e:039})$ and $(\ref{e:040})$ there is a measurable set $G_J\subset J$ satisfying
\begin{equation}\label{zze:024}
    \lambda(G_J)\ge\tfrac34|J|
\end{equation}
such that for every $x\in G_J$ there is an irreducible monic polynomials $P\in\Z[x]$ of degree $n+1$ satisfying $(\ref{e:041})$.
\end{lemma}

\bigskip

\begin{proof}
We will essentially follows the proof of Lemma~\ref{l:04} but replace the construction of $P$ with a different procedure that makes use of the ideas from \cite{Bugeaud-?2}.
Let $G_J=J\smallsetminus E_\infty(J,\delta_1)$, where $\delta_1$ is defined the same way as in the proof of Lemma~\ref{l:04}. Then, we have (\ref{e:029}), which implies (\ref{e:024}). Take any $x\in G_J$. Arguing the same way as in Lemma~\ref{l:04} we obtain $n+1$ linearly independent polynomials $P_j(x)=a_{n,j}x^n+\dots+a_{0,j}\in\Z[x]$ $(0\le j\le n)$ satisfying (\ref{e:032}). The matrix
$A=(a_{i,j})_{0\le i,j\le n}$ satisfies $|\det A| < c_2$ for some constant $c_2$ depending on $n$ only. Again we choose a prime $p$ satisfying (\ref{e:033}) so that $\det
A\not\equiv 0\pmod p$. It is readily verified that
$$
\det(P^{(i)}_j(x))_{0\le i,j\le n}=\det A\prod_{i=0}^n i!\not=0.
$$
Therefore, there is a unique solution $(t_0,\dots,t_n)\in\R^{n+1}$ to the following system of linear equations
\begin{equation}\label{zze:034}
 \frac{(n+1)!}{(n+1-i)!}x^{n+1-i}+p\sum_{j=0}^n t_j P_j^{(i)}(x)=2(n+1)pc_1\xi_i\qquad (0\le i\le n)
\end{equation}
Since $\det A\not\equiv 0\pmod p$ at least one of $a_{0,0},\dots,a_{0,n}$ is not divisible by $p$. Without loss of generality we will assume that $a_{0,0}\not\equiv0\pmod p$. For $j=1,\dots,n$ define $\eta_j=[t_j]$, where $[\,\cdot\,]$ denotes the integer part. Further, define $\eta_0$ to be either $[t_0]$ or $[t_0]+1$ so that
\begin{equation}\label{vb10}
\eta_0a_{0,0}+\dots+\eta_na_{0,n}\not\equiv0\pmod p.
\end{equation}
This is possible because $a_{0,0}\not\equiv0\pmod p$. Define
$$
 P(x)=x^{n+1}+a_nx^n+\dots+a_0:=x^{n+1}+p\sum_{i=0}^n\eta_iP_i(x)\in\Z[x].
$$
Obviously $\deg P=n+1$. The leading coefficient of $P$ is 1 and so is not divisible by $p^2$. By (\ref{vb10}), $a_0\not\equiv0\pmod {p^2}$. However, by the construction of $P$, we have that $a_i\equiv0\pmod p$ for all $i=\overline{1,n}$. Therefore, by Eisenstein's criterion, $P$ is
irreducible over $\Q$.

Finally, it follows from the definition of $\eta_j$ that $|t_j-\eta_j|\le1$ for all $j=0,\dots,n$. Therefore, using the definition of $P$ and (\ref{e:032}) we verify that
$$
    \left|\frac{(n+1)!}{(n+1-i)!}x^{n+1-i}+p\sum_{j=0}^n t_j P_j^{(i)}(x)-P^{(i)}(x)\right|\le(n+1)pc_1\xi_i\qquad (0\le i\le n).
$$
Combining this with (\ref{zze:034}) gives
$$
(n+1)pc_1\xi_i\le |P^{(i)}(x)|  \le  3(n+1)pc_1\xi_i\qquad (0\le i\le n).
$$
Thus, taking $\delta_0=(n+1)pc_1$ and $c_0=3(n+1)pc_1$ gives (\ref{e:041}). The proof is complete.
\end{proof}

\section{Proof of Theorem~\ref{t:02}}\label{proof2}

We now give a complete proof of the theorem in the case $\A^\circ_{n,\nu}(Q,\mu)=\A_{n,\nu}(Q,\mu)$. At the end of the section we will say in what way the proof has to be modified in order to establish the theorem in the case $\A^\circ_{n,\nu}(Q,\mu)=\A^*_{n,\nu}(Q,\mu)$.

\bigskip

Fix $n\ge 2$ and let $\mu$ satisfy (\ref{e:007}). Let $\delta_0$ and $c_0$ be the same as in Lemma~\ref{l:04}. Define the following parameters:
\begin{equation}\label{vb4}
\xi_0=\eta Q^{-n+\mu},\quad \xi_1=\eta^{-n}Q^{1-\mu},\quad \xi_i=\eta Q\quad (2\le i\le n),
\end{equation}
where $0<\eta<1$ is a sufficiently small fixed parameter depending on $n$ only which will be specified later. Fix any interval $J\subset\sui $ and let $\ve=\ve(n,J)$ be the same as in Lemma~\ref{l:04}. Then, (\ref{e:039}) is satisfied with $m\in\{1,2\}$ for sufficiently large $Q$. Also the validity of (\ref{e:040}) easily follows from (\ref{vb4}). Let $G_J$ be the set arising from Lemma~\ref{l:04} and $x\in G_J$. Then, by Lemma~\ref{l:04}, there is a primitive irreducible polynomial $P\in\Z[x]$ of degree $n$ satisfying (\ref{e:041}).

\bigskip

\noindent\textit{Finding $\alpha_1$.}
Let $y\in\R$ be such that $|y-x|=Q^{-n-1+2\mu}$. By (\ref{e:007}), we have $|y-x|<1$. Further, by Taylor's formula,
\begin{equation}\label{vb5}
    P(y)=\sum_{i=0}^n\tfrac{1}{i!}P^{(i)}(x)(y-x)^i.
\end{equation}
Using the inequality $|x-y|<1$, \eqref{e:007}, \eqref{e:041}, and \eqref{vb4} we verify that
\begin{equation}\label{vb+3}
 \Big|P^{(i)}(x)(y-x)^i\Big|< \eta c_0 Q^{-n+\mu}\qquad\text{for }i\ge 2\,.
\end{equation}
Also, by (\ref{e:041}) and (\ref{vb4}), $|P(x)|\le \eta c_0 Q^{-n+\mu}$. Therefore,
\begin{equation}\label{vb6}
 \sum_{i\not=1}\Big|\tfrac{1}{i!}P^{(i)}(x)(y-x)^i\Big|\le
 \eta c_0 Q^{-n+\mu}\sum_{i=0}^n\tfrac{1}{i!}< 3 \eta c_0 Q^{-n+\mu}.
\end{equation}
On the other hand,
\begin{equation}\label{vb7}
|P'(x)(y-x)|\stackrel{\eqref{e:041}\&\eqref{vb4}}{\ge}
\delta_0 \eta^{-n}Q^{-\mu+1} Q^{-n-1+2\mu}\ge\delta_0\eta^{-2} Q^{-n+\mu}.
\end{equation}
It follows from (\ref{vb6}) and (\ref{vb7}) that $P(y)$ has different signs at the endpoints of the interval $|y-x|\le Q^{-n-1+2\mu}$ provided that $\eta\le\delta_0/(3c_0)$. By the continuity of $P$, there is a root $\alpha_1$ of $P$ in this interval, that is
\begin{equation}\label{vb8}
    |x-\alpha_1|< Q^{-n-1+2\mu}\,.
\end{equation}

\noindent\textit{Finding $\alpha_2$.}
Let $y_\rho=x+\rho Q^{-\mu}$, where $2\le|\rho|< Q^{\mu/2}$. In what follows we will again use (\ref{vb5}), this time with $y=y_\rho$. Using $|x-y|<1$, $|\rho|\le Q^{\mu/2}$, \eqref{e:041}, and \eqref{vb4} we verify that
\begin{equation}\label{vb-1}
 \Big|P^{(i)}(x)(y_\rho-x)^i\Big| <  \eta |\rho| c_0Q^{1-2\mu}\qquad\text{for }i\ge 3\,.
\end{equation}
By (\ref{e:007}), (\ref{e:041}), (\ref{vb4}) and the fact that $|\rho|\ge 2$, we have that
$$
|P(x)|\le \eta c_0 Q^{-n+\mu}\le|\rho|\eta c_0 Q^{1-2\mu}
$$
and
$$
|P'(x)(y_\rho-x)|\le \eta^{-n}c_0 Q^{1-\mu}|\rho| Q^{-\mu}=\eta^{-n}c_0|\rho|Q^{1-2\mu}.
$$
The latter two estimates together with (\ref{vb-1}) give
\begin{equation}\label{vb6+}
 \sum_{i\not=2}\Big|\tfrac{1}{i!}P^{(i)}(x)(y_\rho-x)^i\Big|\le
 \eta^{-n} |\rho|c_0 Q^{1-2\mu}\sum_{i=0}^n\tfrac{1}{i!}< 3\eta^{-n} |\rho| c_0 Q^{1-2\mu}.
\end{equation}
On the other hand,
\begin{equation}\label{vb7+}
|\tfrac1{2!}P''(x)(y_\rho-x)^2|\stackrel{\eqref{e:041}\&\eqref{vb4}}{\ge}
\tfrac12\delta_0 \eta Q |\rho|^2 Q^{-2\mu}=\tfrac12\delta_0\eta Q^{1-2\mu}\rho^2.
\end{equation}
It follows from (\ref{vb6+}) and (\ref{vb7+}) that $P(y)$ has the same signs at the points $y_{\pm\rho_0}$ (same as $P''(x)$) with $\rho_0=8c_0\eta^{-n-1}\delta_0^{-1}$.

On the other hand, using (\ref{e:007}) and arguing the same way as during ``Finding $\alpha_1$'', one readily verifies that $P(y_2)$ and $P(y_{-2})$ have different signs. Therefore, $P(y)$ changes sign on one of the intervals
$$
[-\rho_0Q^{-\mu},-2Q^{-\mu}]\qquad\text{or}\qquad [2Q^{-\mu},\rho_0Q^{-\mu}].
$$
By the continuity of $P$, there is a root $\alpha_2$ of $P$ in that interval, that is
\begin{equation}\label{vb8+}
    2Q^{-\mu}\le|x-\alpha_2|< \rho_0Q^{-\mu}\,.
\end{equation}
Combining (\ref{e:007}), (\ref{vb8}) and (\ref{vb8+}) gives $Q^{-\mu}\le |\alpha_1-\alpha_2|\le (\rho_0+1) Q^{-\mu}$, thus establishing (\ref{ne:001}).

\bigskip

\noindent\textit{Estimates for the height.}
Using the fact that $|x|\le\tfrac12$, (\ref{e:041}), (\ref{e:007}) and (\ref{vb4}) we verify that
$$
\begin{array}{l}
 |a_n|\asymp Q,\\[1ex]
 |a_{n-1}|=|P^{(n-1)}(x)-\frac{n!}{1!(n-1)!}\,a_nx|\ll Q,\\[1ex]
 |a_{n-2}|=|P^{(n-2)}(x)-\frac{n!}{2!(n-2)!}\,a_nx^2-\frac{(n-1)!}{1!(n-2)!}\,a_{n-1}x|\ll Q,\\[1ex]
 \vdots
\end{array}
$$
The upshot is that $H(\alpha_1)\asymp Q$. This establishes (\ref{ne:000}) and completes the proof of Theorem~\ref{t:02} in the case $\A^\circ_{n,\nu}(Q,\mu)=\A^*_{n,\nu}(Q,\mu)$.

\bigskip

In the case $\A^\circ_{n,\nu}(Q,\mu)=\A^*_{n,\nu}(Q,\mu)$ the proof remains essentially the same. The only necessary modification arises from taking into account the $(n+1)$-st derivative of $P$. This derivative identically equals $(n+1)!$ and will course no troubles in establishing estimates (\ref{vb6}), (\ref{vb7}), (\ref{vb6+}) and (\ref{vb7+}) which are the key to finding $\alpha_1$ and $\alpha_2$. As to the height, it will be estimated in exactly the same way.

\begin{remark}\label{rem3}
From the above proof we have that $|a_n|\asymp Q$. This condition can be readily used to show that any $\alpha_i$ conjugate to $\alpha_1$ is bounded by a constant depending on $n$ only. This follows from the well known property that $|\alpha_i|\ll H(\alpha_i)/|a_n|$ -- see \cite{Sprindzuk-1969-Mahler-problem}.
\end{remark}

\section{Applications to metric Diophantine approximation}

We begin by recalling a result due to Bernik, Kleinbock and Margulis. In order to state their theorem we introduce the set
$$
\cP_n(\mu,w)=\Big\{x\in[-\tfrac12,\tfrac12]:\left\{\begin{array}{l}
                                                |P(x)|  < H(P)^{-w-\mu}\\[1ex]
                                                |P'(x)| < H(P)^{1-\mu}
                                               \end{array}
\right.\text{ holds for i.m. }P\in\Z[x],\deg P\le n\Big\},
$$
where $n\ge 2$, $\mu\ge0$, $H(P)$ denotes the (absolute) height of $P$ and `i.m.' means `infinitely many'. Applying Dirichlet's pigeonhole principle readily implies that $\cP_n(\mu,w)=[-\tfrac12,\tfrac12]$ if $w\le n-2\mu$. However, when $w> n-2\mu$, the makeup of the set $\cP_n(\mu,w)$ changes completely. The following is a consequence of the Theorem from \cite[\S8.3]{Bernik-Kleinbock-Margulis-01:MR1829381}.

\bigskip

\begin{TheoremB}[(Bernik, Kleinbock and Margulis)]\label{BKM}
Let $n\ge 2$ and $\mu\ge 0$. Then for any $w>n-2\mu$
the set $\cP_n(\mu,w)$ is of Lebesgue measure zero.
\end{TheoremB}

\bigskip

Theorem~BKM is a delicate generalisation of Mahler's problem  \cite{Mahler-1932b} which corresponds to the case $\mu=0$ of Theorem~BKM, though Mahler's problem was settled by Sprind\v zuk \cite{Sprindzuk-1965-Proof-of-Mahler}. It is counterintuitive that for a fixed $\mu$ the set $\cP_n(\mu,w)$ must get smaller as $w$ increases. Hausdorff dimension is traditionally used within this sort of questions in metric Number theory. Using Theorem~\ref{t:02} we are able to produce the following lower bound on the size of $\cP_n(\mu,w)$. In what follows `$\dim$' denotes  Hausdorff dimension.

\bigskip

\begin{theorem}\label{t:03}
Let $n\ge 2$ be an integer and $0<\mu<\frac{n+1}3$. Then
for any $w>n-2\mu$
\begin{equation}\label{e:014}
\dim\cP_n(\mu,w) \ \ge \ \frac{n+1-2\mu}{w+1}.
\end{equation}
\end{theorem}

\bigskip

In the case $\mu=0$, inequality (\ref{e:014}) was first established by Baker and Schmidt \cite{BakerSchmidt-1970} who also conjectured that (\ref{e:014})${}_{\mu=0}$ is actually an equality. This conjecture was proved in \cite{Bernik-1983a}. In view of Theorem~BKM and indeed Theorem~\ref{t:03} it is natural to consider the following generalisation of the Baker-Schmidt conjecture.

\bigskip

\begin{conjecture}\label{c2+}
Let $n$, $\mu$ and $w$ be as in Theorem~\ref{t:03}. Then (\ref{e:014}) is an equality.
\end{conjecture}

\bigskip

Another consequence of this work in the spirit of Theorem~BKM in the following

\begin{theorem}\label{t:04}
Let $v_0,\dots,v_{m-1}\ge0$, $v_m,\dots,v_n\le 0$, $v_0>0$, $v_n<0$
and $v_0+\dots+v_{n}>0$.
Then for almost every $x\in\R$ there are only finitely many $Q\in\N$ such that
\begin{equation}\label{vb9}
         |P^{(i)}(x)| < Q^{-v_i}\quad (0\le i \le n)\qquad \text{for some $P\in\Z[x]\nz$, $\deg P\le n$.}
\end{equation}
\end{theorem}

\bigskip

\begin{proof}
Let $v_0,\dots,v_n$ be given. Without loss of generality we can assume that $x\in\sui$. Let
$$
S_t:=\Big\{x\in\sui:|P^{(i)}(x)| \ll 2^{-v_it}\quad (0\le i \le n)\quad \text{for some $P\in\Z[x]\nz$, $\deg P\le n$}\Big\}.
$$
It is readily seen that our goal is to prove that $\limsup_{t\to\infty}S_t$ has measure zero. By the Borel-Cantelli Lemma, this will follow on showing that $\sum_{t=1}^\infty\lambda(S_t)<\infty$. The latter is easily verified by applying Lemma~\ref{l:03}.
\end{proof}

\bigskip

Mahler's problem corresponds to Theorem~\ref{t:04} with $v_1=\dots=v_{n}=-1$. Theorem~BKM follows from Theorem~\ref{t:04} on taking $v_2=\dots=v_{n}=-1$. Although we are rather flexible in choosing the exponents $v_i$ in Theorem~\ref{t:04}, there are some restrictions which we believe can be safely removed. This is now stated in the form of the following unifying

\bigskip

\begin{conjecture}\label{c2}
Let $\ve>0$. Then for almost every $x\in\R$ the inequality
\begin{equation}\label{vb9+}
         \prod_{i=0}^n|P^{(i)}(x)| < H(P)^{-\ve}
\end{equation}
has only finitely many solutions $P\in\Z[x]$, $\deg P\le n$.
\end{conjecture}

\bigskip

It is likely that in (\ref{vb9+}) the height $H(P)$ can be replaced with $\Pi_+(P):=\prod_{i=1}^n\max\{1,|a_i|\}$, where $P(x)=a_nx^n+\dots+a_1x+a_0$. Also using the inhomogeneous transference principle of \cite{Beresnevich-Velani-08-Inhom} one should be able to establish an inhomogeneous version of Conjecture~\ref{c2} modulo the homogeneous statement.

\subsection{Proof of Theorem~\ref{t:03}}

We will use the ubiquitous systems technique, which is now briefly recalled in a simplified form (see \cite{Beresnevich-Dickinson-Velani-06:MR2184760} for more details and \cite{Beresnevich-Bernik-Dodson-02:MR1975457} for the related notion of regular systems).
Let $I$ be an interval in $\R$ and $\cR:=(\ra)_{\alpha\in J}$ be a family of points $\ra$ in $I$ indexed by a countable set $\cJ$. Let $\beta:\cJ\to \Rp:\alpha\mapsto\ba$ be a function on $\cJ$, which attaches a `weight' $\ba$ to points $\ra$. For $t\in\N$, let $\cJ(t):=\{\alpha \in \cJ:\ba\le 2^t\}$ and assume $\cJ(t)$ is always finite.

Let $\rho: \Rp \to\Rp$ be a function such that
$\lim_{t\to\infty}\rho(t)=0$ referred to as \emph{ubiquity function}. The system $(\cR;\beta)$ is called
{\em locally ubiquitous in $I$ relative to $\rho$} if there is an
absolute constant $k_0>0$ such that for any interval $J\subset I$
\begin{equation}\label{u:00}
\liminf_{t\to\infty}\,\,\lambda\Big(\,\bigcup_{\alpha\in \cJ(t)}
B\big(\ra,\rho(2^t)\big)\cap J\Big) \ \ge \  k_0 \, |J| \, .
\end{equation}

Given a function $\Psi:\Rp\to\Rp$, let
$$
\Lambda_\cR(\Psi) \ := \ \{x\in I:|x-\ra|<\Psi(\ba)
\ \mbox{holds for\ infinitely\ many\ }\alpha \in \cJ \} \,.
$$
The following lemma is Theorem~10 in \cite{Beresnevich-Dickinson-Velani-07:MR2373145}, or alternatively it follows from the more general Corollary~4 from
\cite[p.20]{Beresnevich-Dickinson-Velani-06:MR2184760}.

\begin{lemma}\label{ul}
Let $\Psi:\Rp\to\Rp$ be a monotonic function such that for some
$\phi<1$, $\Psi(2^{t+1})\le \phi\Psi(2^t)$ holds for $t$
sufficiently large. Let $(\cR,\beta)$ be a locally ubiquitous system
in $B_0$ relative to $\rho$. Then for any $s \in (0,1)$
\begin{equation}\label{u:01}
 \cH^s\big(\Lambda_\cR(\Psi)\big) \
= \
\infty \qquad\text{if}\qquad\sum_{t=1}^{\infty}\frac{\Psi(2^t)^s}{\rho(2^t)}\
= \ \infty \, .
\end{equation}
\end{lemma}

\bigskip

\noindent\textit{The ubiquitous system.} Let $n\ge 2$ and $\mu$ satisfy $0<\mu<\frac{n+1}3$. Choose $\mu'=\mu+\delta<\frac{n+1}3$ with $\delta>0$. Let $\cR$ be the set of algebraic numbers $\alpha_1\in\R$ of degree $n$ such that
\begin{equation}\label{ne:001+}
 |\alpha_1-\alpha_2|\le\nu^{-1} H(\alpha_1)^{-\mu'}\qquad\text{for some $\alpha_2\in\R$, conjugate to $\alpha_1$}
\end{equation}
and
\begin{equation}\label{ne:002+}
 |\alpha_i|\ll \nu^{-1}\qquad\text{for any $\alpha_i\in\C$, conjugate to $\alpha_1$},
\end{equation}
where the constant implied by the Vinogradov symbol depends on $n$ only. We will identify $\cJ$ with $\cR$, so that formally $\ra=\alpha$. Further, let $\ba=\nu H(\alpha)$ and $\rho(q):=q^{-n-1+2\mu'}$. Then, by Theorem~\ref{t:02} together with Remark~\ref{rem3}, there is a constant $\nu$ such that $(\cR,\beta)$ is locally ubiquitous in $I:=\sui$ with respect to the above $\rho$. Given $w>0$, let $\Psi(q)=q^{-w-1}$. Clearly, $\Psi(2^{t+1})\le\tfrac12\Psi(2^t)$ and so Lemma~\ref{ul} is applicable to this $\Psi$. Let $s=\frac{n+1-2\mu'}{w+1}$. Since $w>n-2\mu'$, $s<1$. Then
$\frac{\Psi(q)^s}{\rho(q)}$ is identically $1$ and therefore the sum in (\ref{u:01}) diverges. By Lemma~\ref{ul}, we have that
$\cH^s\big(\Lambda_\cR(\Psi)\big)=\infty$.
By the definition of Hausdorff dimension,
$
\dim\Lambda_\cR(\Psi)\ge s=\frac{n+1-2\mu-2\delta}{w+1}.
$
Since $\delta>0$ is arbitrary, it remains to show that
\begin{equation}\label{vbl}
    \Lambda_\cR(\Psi)\subset \cP_n(\mu,w).
\end{equation}
By definition, for every $x\in\Lambda_\cR(\Psi)$ there are infinitely many real algebraic numbers $\alpha_1$ of degree $n$ satisfying (\ref{ne:001+}), (\ref{ne:002+}) and
\begin{equation}\label{ne:003+}
|x-\alpha_1|\ll H(\alpha_1)^{-w-1}.
\end{equation}
Let $P$ denote the minimal polynomial of $\alpha_1$. Then, $P(x)=a_n(x-\alpha_1)\dots(x-\alpha_n)$. By (\ref{ne:001+}), (\ref{ne:002+}), (\ref{ne:003+}) and the fact that $|a_n|\le H(P)$, we get $|P(x)|\ll H\cdot H(P)^{-w-1} H(P)^{-\mu'}= H(P)^{-w-\mu'}$. Since $\mu'>\mu$, we have $|P(x)|< H(P)^{-w-\mu}$ for sufficiently large $H(P)$. Further,
\begin{equation}\label{zzz1}
P'(x)=a_n\sum_{i=1}^{n}\frac{(x-\alpha_1)\dots(x-\alpha_n)}{(x-\alpha_i)}.
\end{equation}
Again, by (\ref{ne:001+}), (\ref{ne:002+}), (\ref{ne:003+}) and the fact that $|a_n|\le H(P)$, we get that every summand in (\ref{zzz1}) is $\ll H(P)^{-\mu'}$, further implying that $|P'(x)|\ll H(P)^{-\mu'}$. Since $\mu'>\mu$, we have $|P'(x)|< H(P)^{-\mu}$ for sufficiently large $H(P)$. The upshot is that the inequalities $|P(x)|< H(P)^{-w-\mu}$ and $|P'(x)|< H(P)^{-\mu}$ hold simultaneously for infinitely many $P\in\Z[x]$ of degree $n$. Thus (\ref{vbl}) is established and the proof is complete.

\section{Final remarks}

\hspace*{\parindent}%
1.
The main body of this paper deals with integral polynomials of degree $n$. However, one can equally develop a similar theory for linear forms of linearly independent analytic functions. This is due to the fact that Lemma~\ref{l:01}, the underlaying fact for other results, is established for linear forms of analytic functions.

\medskip

2.
Clearly, using the algebraic integers part of Theorem~\ref{t:02} it is possible to establish an analogue of Theorem~\ref{t:03} for monic polynomials. Furthermore, using the inhomogeneous transference of \cite{Beresnevich-Velani-08-Inhom} it is possible to establish the inhomogeneous version of Theorem~\ref{t:04}, in particular, the one for monic polynomials.

\medskip

3.
Theorem~\ref{t:02} can be used to give quantitative estimates for the number of polynomials with bounded discriminant -- see, for example, \cite{Bernik-Gotze-Kukso-08:MR2457267}. We are going to address this question in more details in a forthcoming paper.

\medskip

4.
Alongside the Hausdorff dimension generalisation of Theorem~BKM it is interesting to develop a Khintchine type theory -- see \cite[\S8.3]{Bernik-Kleinbock-Margulis-01:MR1829381} where the corresponding problem was stated. When $0<\mu<\frac12$, a result of this kind has been obtained by Kukso \cite{Kukso-07:MR2397717} in the so-called case of divergence.

\medskip

5.
The statement of Theorem~\ref{t:02} can be viewed at a different angle: the algebraic points $(\alpha_1,\alpha_2)$ satisfying (\ref{ne:001}) lie at the distance $Q^{-\mu}$ from the bisector
$y=x$ of the first quadrant. The naturally arising problem is to investigate the distribution of $(\alpha_1,\alpha_2)$ near other rational lines, e.g. $y=2x$ or $y=\frac12x$. More general (and challenging) problem is to investigate the distribution of algebraic points $(\alpha_1,\alpha_2)$ with conjugate coordinates of degree $n$ near non-degenerate curves in the plane, e.g. the parabola $y=x^2$.

\medskip

6.
It would be interesting to develop the theory for non-archimedean extensions of $\Q$ and for `proper' complex algebraic numbers -- see
\cite{Bernik-Gotze-Kukso-08:MR2470800} for a related result.

\medskip

7.
In the previous papers (such as \cite{Bugeaud-Mignotte-2004:MR2096618, Bugeaud-Mignotte-09, Mignotte-83:MR728976}) on the topic of this paper, examples of algebraic numbers with several very close conjugate algebraic numbers have been given. Lemmas~\ref{l:04} and \ref{l:05} of this paper may be used to shed a further light onto this technically more involved question, which will be addressed in a subsequent paper.

\bigskip

\noindent\textit{Acknowledgements.}
The first and second authors are grateful to the University of Bielefeld, where the substantial part of this work has been done, for providing a stimulative research environment during their visits supported by SFB701. The research is also supported by the Royal Society project ``Effective methods in metrical Diophantine approximation''. The authors are also grateful to Yann Bugeaud and the anonymous referee for their very useful comments on an earlier version of this paper.

{\small


\def\cprime{$'$} \def\cprime{$'$} \def\cprime{$'$} \def\cprime{$'$}
  \def\cprime{$'$} \def\cprime{$'$} \def\cprime{$'$} \def\cprime{$'$}
  \def\cprime{$'$} \def\cprime{$'$} \def\cprime{$'$} \def\cprime{$'$}
  \def\cprime{$'$} \def\cprime{$'$} \def\cprime{$'$}
  \def\polhk#1{\setbox0=\hbox{#1}{\ooalign{\hidewidth
  \lower1.5ex\hbox{`}\hidewidth\crcr\unhbox0}}} \def\cprime{$'$}
  \def\cprime{$'$}

}

\

\

\begin{minipage}{\textwidth}
{\small Victor Beresnevich}\\
\footnotesize{\sc University of York, Heslington, York, YO10 5DD,
England}\\
{\it E-mail address}\,: \verb|vb8@york.ac.uk|
\end{minipage}

\vspace*{3ex}

\begin{minipage}{\textwidth}
{\small Vasili Bernik}\\
\footnotesize{\sc Institute of mathematics, Surganova 11, Minsk, 220072, Belarus}\\
{\it E-mail address}\,: \verb|bernik@im.bas-net.by|
\end{minipage}

\vspace*{3ex}

\begin{minipage}{\textwidth}
{\small Friedrich G\"otze}\\
\footnotesize{\sc University of Bielefeld, 33501, Bielefeld, Germany}\\
{\it E-mail address}\,: \verb|goetze@math.uni-bielefeld.de|
\end{minipage}

\end{document}